\documentclass{amsart}
\usepackage{amsfonts}
\usepackage{amsmath,amssymb}
\usepackage{amsthm}
\usepackage{amscd}
\usepackage{graphics}
\usepackage{graphicx}

\theoremstyle{definition}{
\newtheorem{Def}{{\rm Definition}}

\newtheorem{Rem}{{\rm Remark}}

\newtheorem*{MainProb}{Main Problem}
}
\theoremstyle{plain}
{

\newtheorem*{MainThm}{Main Theorem}
\newtheorem*{MainCor}{Main Corollary}

}

\begin{document}
\title[Smooth functions on $3$-dimensional closed manifolds and preimages]{Smooth functions with simple structures on $3$-dimensional closed manifolds with prescribed Reeb graphs and preimages}
\author{Naoki Kitazawa}
\address{Institute of Mathematics for Industry, Kyushu University, 744 Motooka, Nishi-ku Fukuoka 819-0395, Japan\\
 TEL (Office): +81-92-802-4402 \\
 FAX (Office): +81-92-802-4405 \\
}
\email{n-kitazawa@imi.kyushu-u.ac.jp}
\urladdr{https://naokikitazawa.github.io/NaokiKitazawa.html}
\maketitle
\begin{abstract}
We give a new answer to so-called realization problems of graphs as {\it Reeb graphs} of smooth functions with prescribed preimages of regular values having nice structures. We present a best possible answer for functions on $3$-dimensional closed manifolds.

The {\it Reeb space} of a smooth function is the quotient space of the manifold of the domain induced from the following equivalence relation; two points in the manifold are equivalent if and only if they are points of a same connected component of some preimage. They are in considerable cases graphs ({\it Reeb graphs}). 

Reeb spaces with preimages represent the manifolds well and are important tools in geometry. 
Recently they play important roles in applications of mathematics such as visualizations. Realization problems ask us whether we can construct smooth functions with prescribed Reeb graphs and preimages. Studies on construction respecting preimages were essentially started by the author. 
  
\end{abstract}

\section{Introduction}
\label{sec:1}
For a differentiable manifold $X$, $T_xX$ denotes the tangent space at $x$. 
Let $c:X \rightarrow Y$ be a differentiable map. For $x \in X$, ${dc}_x:T_xX \rightarrow T_{c(x)}Y$ denotes the differential at $x$. A point $x \in X$ in the manifold of the domain is said to be a {\it singular} point if ${\rm rank} \ {dc}_x$ is smaller than $\min \{\dim X,\dim Y\}$. 
We define the {\it singular set} of $c$ as the set of all singular points of $c$.
Points in the image of the singular set of $c$ are called {\it singular values} of $c$.
{\it Regular values} of $c$ are points in the manifold of the target which are not singular values.

Differentiable maps are in most cases smooth maps or differentiable maps which we can define the $k$-th differentials for any non-negative integer $k$ at any point in the manifolds of the domains in the present paper. 

${\mathbb{R}}^k$ denotes the $k$-dimensional Euclidean space for $k \geq 1$, endowed with the standard Euclidean metric and for $x \in {\mathbb{R}}^k$, $||x|| \geq 0$ denotes the distance between the origin $0 \in {\mathbb{R}}^k$ and $x$.
$S^k:=\{x \mid x \in {\mathbb{R}}^{k+1}, ||x||=1.\}$ and $D^k:=\{x \mid x \in {\mathbb{R}}^{k}, ||x|| \leq 1.\}$ denote the $k$-dimensional unit sphere and the $k$-dimensional unit disk, respectively.

A {\it Morse-Bott} function $c$ is a smooth function on an $m$-dimensional smooth manifold satisfying the following two.
\begin{itemize}
\item Singular points are always in the interior of the manifold. 
\item At each singular point it has the form represented by the composition of a smooth submersion with a local function represented by the form $$(x_1,\cdots,x_{m^{\prime}}) \mapsto {\Sigma}_{j=1}^{m^{\prime}-i(p)} {x_j}^2-{\Sigma}_{j=1}^{i(p)} {x_{m^{\prime}-i(p)+j}}^2+c(x)$$
for an integer $1 \leq m^{\prime} \leq m$, another suitable integer $0 \leq i(p) \leq m$, and suitable coordinates. 
\end{itemize}
If we do not need a submersion in the composition or equivalently, we can regard this as a diffeomorphism, then the function is a so-called {\it Morse} function.
Singular sets are always closed and smooth submanifolds with no boundaries. If we choose a small open neighborhood of each point which is diffeomorphic to the interior of the unit disk in each connected component of the singular set and orient this and the line of the target around $c(p)$, then $i(p)$ is shown to be taken uniquely in such a way that it is compatible with the local orientations.

A smooth function whose graph is $\{(t,\pm t^2+c)\mid -u<t<u \}$ for a (small) positive number $u>0$ and a real number $c$ is a Morse function. The image is a so-called parabola.   

A {\it fold} map $c$ is a smooth map on a smooth map between smooth manifolds with no boundaries such that at each singular point $x$ of which it has the form represented by the form $$(x_1,\cdots,x_m) \mapsto (x_1,\cdots,x_{n-1},{\Sigma}_{j=1}^{m-n+1-i(p)} {x_{n-1+j}}^2-{\Sigma}_{j=1}^{i(p)} {x_{m-i(p)+j}}^2)$$
for an integer $0 \leq i(p) \leq \frac{m-n+1}{2}$ and suitable coordinates. $i(p)$ is shown to be unique and the singular set is always a closed and smooth submanifold with no boundary, making the restriction of $c$ there a smooth immersion.
If $n=1$ and the manifold of the target is (diffeomorphic to) $\mathbb{R}$, then the function is Morse and if the manifold of the target is a circle, then the map is a so-called {\it circle-valued Morse} function (map).

Canonical projections of unit spheres are shown to be fold maps such that the singular sets are equators and diffeomorphic to unit spheres, that the restrictions to the singular sets are embeddings and that $i(p)=0$ for any singular point $p$. If the manifold of the target is (diffeomorphic to) $\mathbb{R}$ for a canonical projection and we restrict the projection to the preimage of $\{t \mid t \geq 0.\}$ ($\{t \mid t \leq 0.\}$), then we have a Morse function whose singular set consists of exactly one singular point in the interior and $i(p)=0$ (or $i(p)$ and the dimension of the sphere agree). These arguments yield fundamental exercises on fundamental smooth manifolds and Morse functions.  

For systematic theory of Morse functions and the manifolds, see \cite{milnor} for example. \cite{bott} is a pioneering paper on Morse-Bott functions.
For Morse functions and fold maps, \cite{golubitskyguillemin} is a textbook on related theory on singularities of differentiable maps, \cite{thom} and \cite{whitney} are pioneering papers on so-called {\it generic} smooth maps into the plane on smooth manifolds whose dimensions are greater than or equal to $2$, and \cite{saeki0.1} and \cite{saeki0.2} are pioneering works on fold maps and algebraic topological or differential topological properties of the manifolds in the 1990s and 2000s. 

\subsection{Reeb spaces and Reeb graphs}
For any differentiable map $c:X \rightarrow Y$, we can define an equivalence relation ${\sim}_c$ on $X$: $x_1 {\sim}_{c} x_2$ holds if and only if they are in a same connected component of $c^{-1}(y)$ for some point $y \in Y$.

\begin{Def}
The quotient space $W_c:=X/{\sim}_c$ is said to be the {\it Reeb space} of $c$.
\end{Def}
For the Reeb space of $c$, $q_c:X \rightarrow W_c$ denotes the quotient space. We can also define the map $\bar{c}$ uniquely by the relation $c=\bar{c} \circ q_c$. 

In the present paper we consider cases where the Reeb spaces $W_c$ are graphs such that the vertex sets are the sets of all points $p \in W_c$ with the preimages ${q_c}^{-1}(p)$ containing at least one singular point of the maps $c$. 

This is called the {\it Reeb graph} of $c$.
Reeb spaces are shown to be regarded as Reeb graphs for smooth functions on compact manifolds with finitely many singular values in \cite{saeki}.

For Reeb graphs and Reeb spaces, there exist various studies. \cite{reeb} seems to be one of pioneering papers. 
For fold maps, Reeb spaces are shown to be polyhedra whose dimensions are same as those of the manifolds of the targets. See \cite{kobayashisaeki} and \cite{shiota} for example. They have information of topological invariants much in considerable
situations. \cite{saekisuzuoka} and several papers by the author such as \cite{kitazawa0.1}--\cite{kitazawa0.6} show this for fold maps such that preimages of regular values are disjoint unions of copies of spheres.

\subsection{Some notions on graphs topologized canonically.}
A graph is naturally (PL) homeomorphic to a $1$-dimensional polyhedron.
\begin{Def}
An {\it isomorphism} between two graphs $K_1$ and $K_2$ means a (PL) homeomorphism from $K_1$ to $K_2$ mapping the vertex set of $K_1$ onto the vertex set of $K_2$. 
A continuous real-valued function $g$ on a graph $K$ is said to be a {\it good} function if it is injective on each edge.
\end{Def}

In the present paper, we only consider finite and connected graphs with no loops. It immediately follows that a finite graph has a good function if and only if it has no loops as edges.

\subsection{Main Problem and Main Theorem and Corollary.}
The present paper studies the following problem.
\begin{MainProb}
Let a finite and connected graph with at least one edge which may be a multigraph and has no loops be given. Assume also that closed, connected and orientable surface is assigned to each edge. 

Can we construct a smooth function under additional conditions on singularities on a $3$-dimensional closed manifold satisfying the following two. The assumption on surfaces assigned to edges is for prescribed preimages of functions we construct.
\begin{enumerate}
\item The function induces the Reeb graph isomorphic to the graph.
\item Preimages of regular values are as prescribed.
\end{enumerate}
\end{MainProb}
This is a kind of so-called realization problems of graphs as Reeb graphs of smooth functions of suitable classes. \cite{sharko} is a pioneering paper on this respecting (Reeb) graphs only and constructs explicit smooth functions on closed surfaces for such graphs satisfying some suitable conditions. \cite{kitazawa} is regarded as a pioneering paper by the author on this problem and considers conditions on preimages. \cite{saeki} is regarded as a study motivated by this and concerns construction of smooth functions on manifolds whose dimensions are general with prescribed preimages.

Other related papers are also presented in Remark \ref{rem:1} for example.

In the present paper, we first show the following result as our main theorem.
\begin{MainThm}
\label{mthm}
Let $F_1$ and $F_2$ be non-empty spaces and closed surfaces which may not be connected and the difference of the Euler numbers are even. Let $|F_i|$ denote the number of all connected components of $F_i$ for $i=1,2$.
Let $[a_1,a_2] \subset \mathbb{R}$ be a closed interval.
Then there exist a $3$-dimensional compact and connected manifold ${\tilde{M}}_{F_1,F_2}$ and a smooth function ${\tilde{f}}_{F_1,F_2}$ on ${\tilde{M}}_{F_1,F_2}$ satisfying the following properties.
\begin{enumerate}
\item The image of ${\tilde{f}}_{F_1,F_2}$ is $[a_1,a_2]$.
\item The boundary of ${\tilde{M}}_{F_1,F_2}$ is diffeomorphic to the disjoint union of $F_1$ and $F_2$ and coincides with the preimage of $\{a_1,a_2\}$ and the preimage of $a_i$ is diffeomorphic to $F_i$ for $i=1,2$.
\item The Reeb space $W_{{\tilde{f}}_{F_1,F_2}}$ is PL homeomorphic to a graph whose vertex set consists of exactly $|F_1|+|F_2|+1$ vertices and the number of vertices whose degrees are greater than $1$ is $1$.
\item There exists exactly one singular value $a \in (a_1,a_2)$ for the function ${\tilde{f}}_{F_1,F_2}$.
\item The function is, around each singular point, represented as at least one of the following functions for suitable coordinates.
\begin{enumerate}
\item A Morse function.
\item The composition of a fold map into a surface with a Morse function.
\end{enumerate}
\end{enumerate}
\end{MainThm}
By respecting the preprint \cite{kitazawa3} with \cite{kitazawa}, we also have the following corollary. 
\begin{MainCor}
\label{mcor}
Let $K$ be a finite and connected graph which has at least one edge and no loops. Let $g$ be a good function on $K$. Let $r_K$ be an integer-valued map on the edge set $E$. Suppose that these two functions satisfy the following conditions.
\begin{itemize}
\item For each vertex $v$ where $g$ has a local extremum, the number of edges containing $v$ at which the values of $r_K$ are odd and negative is even.
\item For each vertex $v$ where $g$ does not have a local extremum, let $M_v$ {\rm (}$m_v${\rm )} denote the number of edges containing $v$ such that the values of $r_K$ there are odd and negative and that the restrictions of $g$ there have the maxima {\rm (}resp. minima{\rm )} at $v$. Then $M_v-m_v$ is even.
\end{itemize}
Then there exist a $3$-dimensional closed and connected manifold $M$ and a smooth function $f$ on $M$ satisfying the following properties.
\begin{enumerate}
\item The Reeb graph $W_f$ of $f$ is isomorphic to $K$ and we can take a suitable isomorphism $\phi:W_f \rightarrow K$ compatible with the remaining properties.
\item For each point $\phi(p) \in K$ {\rm (}$p \in W_f${\rm )} in the interior of an edge $e$, the preimage ${q_f}^{-1}(p)$ is a closed, connected, and orientable surface of genus $r_K(e)$ if $r_k(e) \geq 0$ and a closed, connected, and non-orientable surface of genus $-r_K(e)$ if $r_k(e)<0$.
\item For each point $p \in M$ mapped by $q_f$ to a vertex $v_p:=q_f(p) \in W_f$, $f(p)=g \circ \phi(v_p)$.
\item The function is, around each singular point, represented as at least one of the following functions for suitable coordinates.
\begin{enumerate}
\item A Morse function.
\item A Morse-Bott function which is not Morse.
\item The composition of a fold map into a surface with a Morse function.
\item The composition of a Morse function with a Morse function whose graph is a parabola.
\item The composition of a fold map into a surface with a Morse function and a Morse function whose graph is a parabola.
\end{enumerate}
\end{enumerate}
\end{MainCor}
We have shown this in the case where the preimages of regular values are orientable with the fact that we can have the $3$-dimensional manifold as an orientable manifold as the main theorem or Theorem 1 of \cite{kitazawa} for example. 

We prove Main Theorem in the next section. Methods we use are similar to the methods used in the main theorem or Theorem 1 of \cite{kitazawa} and some theorems in \cite{kitazawa2}--\cite{kitazawa4}.
We first construct local smooth maps whose Reeb spaces are (PL) homeomorphic to small regular neighborhoods of vertices and glue the functions together.
A new ingredient of the present study is the construction in STEP C in the proof and iterations of known technique in \cite{kitazawa3} in STEP D. 
\section{Proofs of Main Theorem and Corollary.}
\begin{proof}[A proof of Main Theorem.]
We first review existing arguments and show partially and generalize situations to complete the proof. \\
\ \\
STEP A Reviewing ingredients of proofs of important theorems in \cite{kitazawa} and \cite{kitazawa3} and proving the theorem in the case where $F_1$ and $F_2$ consist of only closed and connected surfaces whose Euler numbers are
 even. \\
 
By virtue of \cite{michalak}, we can show the case where only $0$ is assigned to each edge. We have a desired Morse function whose image is $[a_1,a_2]$. We can remove finitely many disjoint trivial smooth bundles over $[a_1,a_2]$ whose fibers are diffeomorphic to $D^2$ and which are apart from the singular set by virtue of the structure of the function. We attach copies of the Morse function satisfying the following three with a trivial smooth bundle over $[a_1,a_2]$ whose fibers is diffeomorphic to $D^2$ and which is apart from the singular set removed instead, preserving the values, by diffeomorphisms. We can also remove the trivial bundles by similar reasons.   
\begin{itemize}
\item The preimage of $a_1$ or $a_2$ is diffeomorphic to $S^2$.
\item The preimage of $a_1$ or $a_2$ is diffeomorphic to $S^1 \times S^1$ or the Klein Bottle.
\item The singular value set is $\{a\} \subset (a_1,a_2)$.
\end{itemize} 
Choosing a suitable Morse function such that preimages of regular values are disjoint unions of copies of $S^2$, removing the trivial bundles in a suitable way and attaching the copies of the functions, yields a desired Morse function. We can also do the latter two operations by the structures of the smooth functions. \\
\ \\
STEP B Reviewing ingredients of proofs of important theorems in \cite{kitazawa3} and proving the theorem in the case where $F_1$ and $F_2$ are connected surfaces whose Euler numbers are odd. \\

First we consider the case where $F_1$ and $F_2$ are both diffeomorphic to the real projective plane. 
See \cite{milnor} for ({\it $k$-}){\it handles} for a non-negative integer $k$ and relationship between a singular point $p$ and an $i(p)$-handle for Morse functions where $i(p)$ is as in the definition of a Morse function in the introduction. We omit rigorous expositions in the present paper. See \cite{milnor}, introduced also in the introduction.

By attaching two handles to $F_1 \times \{0\} \subset F_1 \times [-1,0]$ and we have a $3$-dimensional compact, connected and smooth manifold whose boundary is diffeomorphic to the disjoint union $F_1 \sqcup F_2$. More precisely, we attach a $1$-handle, diffeomorphic to $D^1 \times D^2 \supset \partial D^1 \times D^2$, to a smooth submanifold $S_1$ diffeomorphic to the disjoint union of two copies of $D^2$ in $F_1 \times \{0\}$ and a $2$-handle, diffeomorphic to $D^2 \times D^1 \supset \partial D^2 \times D^1$, to a smooth submanifold $S_2$ diffeomorphic to $S^1 \times D^1$ in $F_1 \times \{0\}$, apart from the previous submanifold $S_1$ and locating the two copies of the unit disk which are connected components of the submanifold $S_1$ before in different connected components of the complementary set of $S_2$ in $F_1 \times \{0\}$. This yields a desired Morse function with exactly two singular points and exactly one singular value $a$.

For a general case here, we (can) do operations as in STEP A to complete the present step. \\
\ \\
STEP C Proving the theorem in the following case: either $F_1$ or $F_2$ is a disjoint union of exactly two closed and connected surfaces whose Euler numbers are odd and the other surface is a closed, connected and orientable surface or non-orientable one with even Euler number. This is a new ingredient in our proof. \\

First we show this for the case where $F_1$ is diffeomorphic to $S^2$ and $F_2$ is diffeomorphic to the disjoint union of two copies of the projective plane.
We consider a compact and connected surface $S_{a,r}$ in ${\mathbb{R}}^3$ defined in the following way.
\begin{enumerate}
\item We choose an arbitrary sufficiently large positive number $r>0$.
\item We consider three circles $C_1$, $C_2$, and $C_3$ in ${\mathbb{R}}^2 \times \{0\} \subset {\mathbb{R}}^3$ centered at $(0,a+r,0)$ whose radii are $r$, $2r$ and $3r$ respectively. 
\item For each point $p$ in $C_2$ whose second component is smaller than or equal to $a+r$, consider the circle whose radius is $r$, which is centered there, and which passes through the following two points uniquely found in the segment of length $3r$ starting from $(0,a+r,0)$ before and containing $p$.
\begin{enumerate}
\item The point in $C_1$.
\item The point in $C_3$.
\end{enumerate}
We can also find the circle uniquely.
\item We take the disjoint union of all the circles before and we have the desired compact and connected surface $S_{a,r}$ smoothly embedded in ${\mathbb{R}}^3$.
\end{enumerate}
Let $S_{a,a_1,a_2,r} \subset S_{a,r}$ be the subspace of all points whose second components $t_0$ satisfy $a_1 \leq t_0 \leq a_2$. 
We consider the following two hyperplanes in ${\mathbb{R}}^3$ where $\delta>0$ is a positive real number.
\begin{itemize}
\item $H_1:=\{(x_1,x_2,x_3) \mid x_1=-3r+\delta, (x_2,x_3) \in {\mathbb{R}}^2. \}$.
\item $H_2:=\{(x_1,x_2,x_3) \mid x_1=3r-\delta, (x_2,x_3) \in {\mathbb{R}}^2. \}$.
\end{itemize}
Furthermore, we can take $\delta$ so that for $H_j \bigcap S_{a,a_1,a_2,r}$ and each $j=1,2$, the minimum of the second component is $a$, the maximum of this is $a_2$ with exactly two points at which the values of the second components are this maximum. We choose $\delta$ in this way. These subspaces are parabolas and diffeomorphic to a closed interval. $H_{S,j}$ denotes the subspace for $j=1,2$.
We define another parabola in the image of the canonical projection of $S_{a,a_1,a_2,r} \subset {\mathbb{R}}^3$ to ${\mathbb{R}}^2 \times \{0\}$ defined by $(x_1,x_2,x_3) \mapsto (x_1,x_2,0)$. $r$ is sufficiently large. We can define the parabola $\{(t,-(t-u_0)^2+a,0) \mid -u_1 \leq t-u_0 \leq u_1.\}$ by taking suitable real numbers $u_0$ and $u_1>0$ so that $-{u_1}^2+a=a_1$.
We can define a smooth curve in $S_{a,a_1,a_2,r}$ mapped onto the parabola by a diffeomorphism by the canonical projection.
Let $P_S$ denote the curve. We define the set of all points in $C_1$ the values of whose second components are smaller than or equal to $a_2$ and let $C_S \subset S_{a,a_1,a_2,r}$ denote the preimage of the set for the canonical projection. $C_S$ is also a subset of ${\mathbb{R}}^2 \times \{0\}$ and the restriction of the projection to $C_S$ is regarded as the identity map.

Hereafter, let ${\pi}_S$ denote the canonical projection of the surface $S_{a,a_1,a_2,r}$.

We also have a smooth map on a $3$-dimensional compact and connected manifold ${\tilde{M}}_{F_1,F_2}$ into $S_{a,a_1,a_2,r}$ satisfying the following properties.
\begin{itemize}
\item The restriction to the singular set is a smooth embedding and the image of the singular set is $H_{S,1} \sqcup H_{S,2}  \sqcup P_S \sqcup C_S$ .  
\item The restriction to the preimage of the interior of $S_{a,a_1,a_2,r}$ is a fold map such that the restriction to the singular set is an embedding where the surface of the target is restricted to the interior.
\item The restriction to the preimage of each connected component of the boundary of the surface is a circle-valued Morse map where the manifold of the target is restricted to the circle.
\item The image of the map is the complementary set of the disjoint union of the following subspaces.
\begin{itemize}
\item The connected component of the preimage of the subspace $\{(t,t^{\prime}) \mid -u_1 < t-u_0 < u_1, a_1 \leq t^{\prime}<-(t-u_0)^2+a.\}$ for the canonical projection ${\pi}_S$ of the surface of the target of the smooth map whose closure contains $P_S$ in the surface.
\item The preimage of the subspace $\{(x_1,x_2,x_3) \mid x_1<-3r+\delta, (x_1,x_2,x_3) \in {\pi}_S(S_{a,a_1,a_2,r}).\}$ for the canonical projection ${\pi}_S$ of the surface of the target of the smooth map.
\item The preimage of the subspace $\{(x_1,x_2,x_3) \mid x_1>3r-\delta, (x_1,x_2,x_3) \in {\pi}_S(S_{a,a_1,a_2,r}).\}$ for the canonical projection ${\pi}_S$ of the surface of the target of the smooth map.
\end{itemize}
\item Preimages of regular values are always diffeomorphic to a circle or the empty set. 
\end{itemize}
FIGURE \ref{fig:1} shows the composition of the map with the canonical projection to ${\mathbb{R}}^2 \times \{0\} \subset {\mathbb{R}}^3$  
\begin{figure}
\includegraphics[width=60mm,height=60mm]{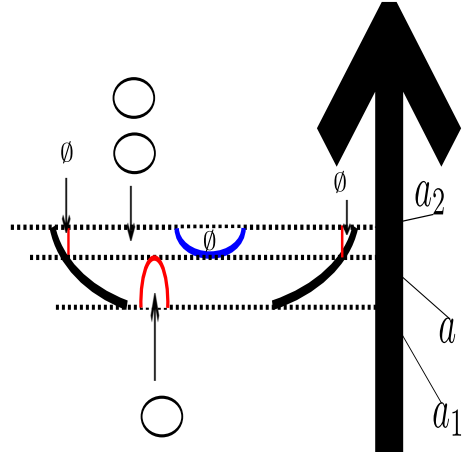}
\caption{The image of the composition of the smooth map into $S_{a,a_1,a_2,r}$ with the canonical projection to ${\mathbb{R}}^2 \times \{0\} \subset {\mathbb{R}}^3$. The red segments are for the images of $H_{S,j}$, the red parabola is for the image of $P_S$, and the blue arc is for the image of $C_S$. Circles and $\emptyset$ are for types of preimages.}
\label{fig:1}
\end{figure}

By composing the map with the canonical projection to the second component yields a desired function. The restriction of the resulting function to the interior is represented as the composition of a fold map into a surface with a Morse function. For the topologies around the preimages, consult fundamental theory on Morse functions on closed surfaces and \cite{saeki0.2}. This is presented again in the proof of Main Corollary. \\ 
We can exchange roles of $F_1$ and $F_2$ of course. \\
For a general case, the structures of the smooth functions enable us to show as in STEP A and STEP B. \\
\ \\
\noindent STEP D Proving the theorem in a general case to complete the proof. \\
\ \\
We construct suitable families of functions as in STEP A, STEP B and STEP C. To obtain a desired smooth function starting from these functions, we apply operations presented in FIGURE \ref{fig:2} and FIGURE \ref{fig:3} in a suitable way. 

\begin{figure}
\includegraphics[width=60mm,height=60mm]{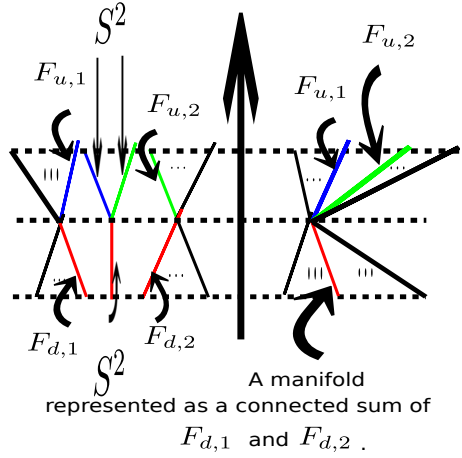}
\caption{Obtaining a new smooth function and its Reeb space by gluing two given functions and another Morse function with suitable trivial smooth bundles over the images whose fibers are diffeomorphic to $D^2$ removed: the preimages of the extrema are diffeomorphic to $S^2$ and $S^2 \sqcup S^2$ respectively for the original Morse function for the last function. The $1$-dimensional objects consisting of straight segments show Reeb spaces of smooth functions and from functions in the left we have the function in the right. $F_{u,1}$, $F_{d,1}$ and $S^2$ show types of preimages for the original functions before we remove the trivial bundles whose fibers are diffeomorphic to $D^2$.}
\label{fig:2}
\end{figure}
\begin{figure}
\includegraphics[width=60mm,height=60mm]{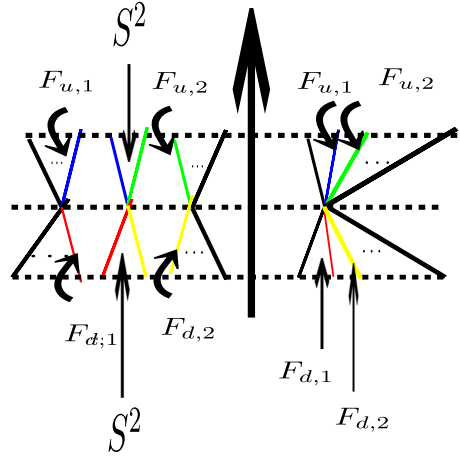}
\caption{Obtaining a new smooth function and its Reeb space by gluing two given functions and another Morse function with suitable trivial smooth bundles over the images whose fibers are diffeomorphic to $D^2$ removed: the preimages of the extrema are both diffeomorphic to $S^2 \sqcup S^2$ for the original Morse function for the last function. $F_{u,1}$, $F_{d,1}$ and $S^2$ show types of preimages for the original functions before we remove the trivial bundles whose fibers are diffeomorphic to $D^2$.}
\label{fig:3}
\end{figure}

They show some fundamental operations deforming Reeb spaces and preimages where descriptions of the types of the manifolds are for preimages for original functions. We remove trivial smooth bundles over the images whose fibers are diffeomorphic to $D^2$ from the original functions apart from the singular sets and glue the resulting functions preserving the value at each point in the manifolds of the domains. 

We first prepare a suitable family of local functions in STEP A, STEP B and STEP C. By a suitable iteration of these two operations, we have a desired function. 

We explain FIGURE \ref{fig:2} and FIGURE \ref{fig:3} shortly. In each situation, let two functions on compact and connected manifolds satisfying the desired properties be given.
FIGURE \ref{fig:2} constructs a new function satisfying the desired properties and the following properties. 
\begin{itemize}
\item Either the preimage of the minimum or the maximum is the disjoint union of the original two preimages of the minima or the maxima.
\item Either the preimage of the minimum or the maximum is the disjoint union of the following two surfaces.
\begin{itemize}
\item The disjoint union of the original two preimages of the minima or the maxima with exactly one connected component removed for each of the given two functions.
\item A surface represented as a connected sum of the removed two surfaces before.
\end{itemize}
\end{itemize}
FIGURE \ref{fig:3} constructs a new function satisfying the desired properties and the following properties.
\begin{itemize}
\item The preimage of the minimum is the disjoint union of the original two preimages of the minima.
\item The preimage of the maximum is the disjoint union of the original two preimages of the maxima.
\end{itemize}
For the figures, see also Lemmas in \cite{kitazawa3} for a more precise exposition. 

We show more detailed arguments to obtain a desired local function.

\begin{itemize}
\item The case where the number of connected components whose Euler numbers are odd of $F_1$ and that of $F_2$ agree.
\begin{itemize} 
\item The case where there exist no connected components whose Euler numbers are odd. \\
We have a desired function as a function in STEP A.
\item The case where there exists at least one connected component whose Euler number is odd. \\
We need at least one function in STEP B.
\begin{itemize}
\item The case where there exist no connected components whose Euler numbers are even. \\
We prepare a suitable family of functions in STEP B. We
apply operations of FIGURE \ref{fig:3} one after another in a suitable way. 
\item The case where there exists at least one connected component whose Euler number is even for both $F_1$ and $F_2$. \\
We prepare a suitable family of functions in STEP B. We
apply operations of FIGURE \ref{fig:3} one after another in a suitable way. We prepare a suitable function in STEP A and apply an operation of FIGURE \ref{fig:3} in a suitable way.
\item The case where there exists at least one connected component whose Euler number is even for either $F_1$ or $F_2$. \\
We prepare a suitable family of functions in STEP B. We
apply operations of FIGURE \ref{fig:3} one after another in a suitable way to have a smooth function. We prepare a suitable function in STEP A such that at least one of the preimage of the minimum and that of the maximum is diffeomorphic to $S^2$. We apply an operation of FIGURE \ref{fig:2} to these two functions in a suitable way.
\end{itemize}
\end{itemize}
\item The case where the number of connected components whose Euler numbers are odd of $F_1$ and that of $F_2$ do not agree.
\begin{itemize}
\item The case where there exist no connected components whose Euler numbers are even. \\
We need at least one function in STEP C. We prepare a suitable family of functions in STEP C. We also pose the condition on each function of these functions that the preimage of the minimum or the maximum is diffeomorphic to $S^2$. 
We apply operations of FIGURE \ref{fig:2} one after another to obtain a new function on a compact and connected manifold such that the preimage of the minimum or the maximum is diffeomorphic to $S^2$. We need at least one function in STEP B, prepare suitable families of functions in STEP B and apply operations of FIGURE \ref{fig:3} one after another to obtain a new function on a compact and connected manifold in a suitable way. Thus we have two functions and for these functions we apply an operation of FIGURE \ref{fig:2} in a suitable way.
\item The case where there exists at least one connected component whose Euler number is even for both $F_1$ and $F_2$.
\begin{itemize}
\item The case where there exists at least one connected component whose Euler number is odd for both $F_1$ and $F_2$. \\
We construct a suitable function on a compact and connected manifold such that the preimages of the minimum and the maximum have no connected components whose Euler numbers are even by the method just before. We prepare a suitable function in STEP A and apply an operation of FIGURE \ref{fig:3} in a suitable way.
\item The case where there exists at least one connected component whose Euler number is odd for either $F_1$ or $F_2$. \\
We need at least one function in STEP C, prepare a suitable family of functions in STEP C, and pose the condition on each function of these functions that the preimage of the minimum or the maximum is diffeomorphic to $S^2$. We apply operations of FIGURE \ref{fig:2} one after another to obtain a new function on a compact and connected manifold such that the preimage of the minimum or the maximum is diffeomorphic to $S^2$. 
We also prepare a suitable function in STEP A such that at least one of the preimage of the minimum or the maximum is connected and diffeomorphic to $F_0$. 
Thus we have two smooth functions. We apply an operation of FIGURE \ref{fig:2} in a suitable way to obtain a new smooth function such that at least one of the preimage of the minimum or the maximum is connected and diffeomorphic to $F_0$. This completes an exposition for the case where the desired preimage of the minimum or the maximum is connected and has the even Euler number. In general, we also need to prepare a suitable function in STEP A such that at least one of the preimage of the minimum or the maximum is diffeomorphic to $S^2$ and apply an operation of FIGURE \ref{fig:2} to the previously obtained function and this function in a suitable way.
\end{itemize}
\item The case where there exists at least one connected component whose Euler number is even for either $F_1$ or $F_2$. \\
\begin{itemize}
\item The case where there exists at least one connected component whose Euler number is odd for both $F_1$ and $F_2$. \\
We construct a suitable function on a compact and connected manifold such that the preimages of the minimum and the maximum have no connected components whose Euler numbers are even as in the first case here.
We prepare a suitable function in STEP A such that the preimage of the minimum or the maximum is diffeomorphic to $S^2$. We
apply an operation of FIGURE \ref{fig:2} to these two functions in a suitable way.
\item The case where there exists at least one connected component whose Euler number is odd for either $F_1$ or $F_2$. \\
We need at least one function in STEP C, prepare a suitable family of functions in STEP C, and
apply operations of FIGURE \ref{fig:2} one after another to obtain a new function on a compact and connected manifold such that the preimage of the minimum or the maximum is diffeomorphic to $S^2$ or another closed and connected surface whose Euler number is even. This completes an exposition for the case where the desired preimage of the minimum or the maximum is connected and has the even Euler number. In general, we need to prepare a suitable function in STEP A such that the preimage of the minimum or the maximum is diffeomorphic to $S^2$ and apply an operation of FIGURE \ref{fig:2} to these two functions in a suitable way.
\end{itemize}
\end{itemize}
\end{itemize}

\end{proof}
We present a proof of Main Corollary. 
A {\it linear} bundle is a smooth bundle whose fiber is diffeomorphic to a Euclidean space, a unit sphere, or a unit disk, and whose structure group acts on the fiber by linear transformations.
\begin{proof}[A proof of Main Corollary.]
As done in important theorems in \cite{kitazawa}--\cite{kitazawa4}, we construct local smooth functions and glue them together by suitable diffeomorphisms between  connected components of the boundaries of the manifolds.\\
\ \\
STEP A Around a vertex $v_p$ of degree 1 at which the given good function $g$ has a local extremum. \\
We need important arguments and theorems of \cite{kitazawa} and \cite{kitazawa3}. 

If the value of $r_K$ at the edge containing the vertex is $-2$, then we consider a Morse Bott-function on a manifold diffeomorphic to the total space of a linear bundle over $S^1$ whose fiber is diffeomorphic to $D^2$ and which is not trivial such that the singular set is $S^1 \times \{0\}$ and mapped to $g(v_p)$. Here $S^1 \times \{0\}$ denotes the image of the section taking values corresponding to the origin $0$ in the unit disk $D^2$. This function is not presented in these two papers.

If the value of $r_K$ at the edge containing the vertex is an even number $-2(k+1)$ ($k>0$), then we have a local map represented as the composition of smooth surjection onto the $2$-dimensional unit disk $D^2$ such that the restriction to the interior is a fold map onto the interior of the disk with a Morse function on the $2$-dimensional disk represented by the form $(x_1,x_2) \mapsto \pm({x_1}^2+{x_2}^2)+g(v_p)$ for suitable coordinates. We show the local map onto $D^2 \subset {\mathbb{R}}^2$ applying the exposition in the proof of Theorem 1 of \cite{kitazawa}.

Let $b>1$. Let $\{s_j\}_{j=1}^k$ be an increasing sequence of real numbers satisfying the following conditions.
\begin{itemize}
\item $s_1s_k \leq 0$. 
\item If $k=1$, then $s_1=0$.
\item If $k>1$, then $s_1<0$ or $s_k>0$.
\end{itemize}
For a compact and connected surface $P$ obtained by removing the interiors of disjointly and smoothly embedded two copies of the $2$-dimensional unit disk $D^2$ in a suitable $2$-dimensional closed and connected surface, we have a smooth map $F:P \times [-b,b] \rightarrow [-b,b] \times [-b,b]$ satisfying the following properties. 

\begin{itemize}
\item We can define a family $\{f_t(u):=F(u,t)\}_{t \in [-b,b]}$ of smooth functions on $P$ each of which is $f_t:[-b,b] \rightarrow [-b,b]$.
\item The images of these functions are $[-b,b]$ and these functions have exactly $k$ singular points.
\item On the interior they are Morse functions with exactly $k$ singular points whose indices are $1$ and the singular value set of $f_t$ is $\{s_j-\frac{s_j}{b}(t+b) \mid 1 \leq j \leq k\}$.
\item For each $f_t$, the preimage of each regular value $p \in [-b,b]$ is a circle.
\item The preimage of each singular value of the function $f_t$ is as follows. 
\begin{itemize} 
\item A $1$-dimensional polyhedron obtained by identifying two distinct points in a circle for $f_t$ ($t \neq 0$). 
\item A $1$-dimensional polyhedron obtained in the following way for $f_t$ ($t=0$).
\begin{itemize}
\item Consider a circle and $2k$ distinct points there.
\item Divide the $2k$ distinct points into suitable $k$ pairs of the points.
\item Identify the two distinct points in each of the $k$ pairs on the circle.
\end{itemize}
\end{itemize}
\item (The restriction of) $F$ (to the interior) is a fold map and the singular set of $F$ is a disjoint union of $k$ copies of a closed interval. The restriction to each singular set is an embedding and its image is $\{(s_j-\frac{s_j}{a}(t+a),t) \in [-a,a] \times [-a,a] \mid -a \leq t \leq a, 1 \leq j \leq k\}$.
\end{itemize}   

We restrict the map to the preimage of the $D^2 \subset {\mathbb{R}}^2$.

Consult \cite{kitazawa} and \cite{kitazawa3} for more precise expositions. For theory of {\it fibers}, see \cite{saeki0.2}. This concerns theory of germs of smooth maps around preimages of points.\\
\ \\
STEP B Around a vertex at which the given good function does not have a local extremum. \\
We have this immediately from Main Theorem. 
Note that the Euler number of a closed and connected surface is odd if and only if it is non-orientable and has an odd genus. \\
\ \\
STEP C Around a vertex $v_p$ of degree 2 at which the given good function $g$ has a local extremum. \\
We construct a function by the proof of Main Theorem by defining $F_1$ and $F_2$ there suitably as non-empty sets and closed surfaces so that the disjoint union is diffeomorphic to a closed surface obtained by considering the disjoint union of closed, connected and orientable surfaces of genera $r_k(e)$ (for $r_K(e) \geq 0$) or closed, connected and non-orientable surfaces of genera $-r_k(e)$ (for $r_K(e)<0$) for all edges $e$ containing $v_p$. 
We need to and can choose these two surfaces so that the numbers of connected components whose Euler numbers are odd are same. We can do this by virtue of the assumption on $g$ and $r_K$. We embed the image into the plane so that the image is a parabola and of the form $\{(t,\pm t^2+g(v_p))\mid -u<t<u\} \subset {\mathbb{R}}^2$ for a small positive real number $u>0$ and that the singular value is mapped to $(0,g(v_p))$. After that we compose the resulting map with the canonical projection to the second component. For similar methods, consult \cite{kitazawa}, \cite{kitazawa2}, \cite{kitazawa3} and \cite{kitazawa4} for example. \\
\ \\
We construct trivial smooth bundles on the complementary set of the interior of the disjoint union of suitable small regular neighborhoods of all vertices. 
\ \\
Last, by gluing the local functions together by suitable diffeomorphisms between boundaries, we have a desired smooth function. We omit explicit expositions on construction of a desired PL homeomorphism $\phi:W_f \rightarrow K$. See \cite{kitazawa}--\cite{kitazawa4}. This completes the proof.
\end{proof}
\begin{Rem}
\label{rem:1}
\cite{saeki} is a study motivated by \cite{kitazawa}. This generalizes compact (closed), connected and smooth manifolds assigned to edges.
We must obey so-called cobordism relations on compact manifolds in assigning compact (closed), connected and smooth manifolds to edges modulo existence of diffeomorphisms between smooth manifolds. For example, in Main Theorem the difference of the Euler numbers of $F_1$ and $F_2$ are assumed to be even and we must not drop the condition. \cite{saeki} presents construction of smooth functions using smooth functions which are not real analytic under the condition that the dimensions of the compact manifolds are general. As \cite{kitazawa}, \cite{kitazawa3} and \cite{kitazawa4} say, in \cite{saeki} explicit singularities of smooth functions obtained there are not studied.

\cite{kitazawa3} (partially) generalizes \cite{kitazawa} and locally it is regarded as a specific study on Main Theorem. Different from our present study, a sufficient condition to obtain a smooth function which is always locally Morse around each singular point is studied. \cite{kitazawa4} is also a study closely related to \cite{kitazawa3}.
\end{Rem}
\section{Acknowledgement, grants and data.}
\label{sec:3}
\thanks{The author is a member of the project JSPS KAKENHI Grant Number JP17H06128 "Innovative research of geometric topology and singularities of differentiable mappings" (Principal Investigator: Osamu Saeki). This work is supported by the project. 
The author would like to thank Osamu Saeki again much for private discussions on \cite{saeki} respecting \cite{kitazawa}. These discussions are encouraging the author to study related studies further including the present study. 
All data essentially supporting the present study are all in the present paper.}

\end{document}